\documentclass[12pt]{amsart}
\usepackage{amscd,amssymb,epsfig}
\usepackage{graphics}

\input epsf
\newtheorem{theorem}{Theorem}[section]
\newtheorem{lemma}[theorem]{Lemma}
\newtheorem{proposition}[theorem]{Proposition}
\newtheorem{corollary}[theorem]{Corollary}
\theoremstyle{definition}

\theoremstyle{remark}
\newtheorem{remark}[theorem]{Remark}
\numberwithin{figure}{section}
\numberwithin{table}{section}

\begin{document} 

\title[Measured Laminations and Convexity] 
{Tropicalized Lambda Lengths, Measured Laminations and Convexity}

\author{R.\ C.\ Penner}
\address{Center for the Quantum Geometry of Moduli Spaces\\
Aarhus University\\
DK-8000 Aarhus C, Denmark\\~~{\rm and}~
Departments of Mathematics and Theoretical Physics\\
Caltech\\
Pasadena, CA 91125\\
USA\\}
\email{rpenner{\char'100}imf.au.dk}

\thanks{This work has been supported by the Centre 
for Quantum Geometry of Moduli Spaces
funded by the Danish National Research Foundation.}
\thanks{The author is grateful to Volodya Fock for stimulating discussions and
especially to Dylan Thurston for his participation and contributions to this paper.}

\begin{abstract}
This work uncovers the tropical analogue for measured laminations of the convex hull construction of decorated 
Teichm\"uller theory, namely, it is a study in coordinates of geometric degeneration to a point of Thurston's boundary for Teichm\"uller space.
This may offer a paradigm for the extension of the basic cell decomposition of  Riemann's moduli space
to other contexts for general moduli spaces of flat connections on a surface.  In any case, this discussion drastically simplifies aspects of previous related studies as is explained.  Furthermore, a new class of measured laminations relative to an ideal cell decomposition of a surface is discovered in the limit.   Finally, the tropical analogue of the convex hull construction in Minkowski space is formulated as an explicit algorithm that serially simplifies a triangulation with respect to a fixed lamination and has its own independent applications.
\end{abstract}

\maketitle

\section*{Introduction}

As we shall recall from \cite{Penner87,Penner11} in detail in the next section,  lambda lengths $\lambda:\Delta\to{\mathbb R}_{>0}$ on the edges $\Delta$ of an ideal triangulation of a
punctured surface $F$ give global affine coordinates on the space $\widetilde {\mathcal T}(F)$ of decorated hyperbolic structures on $F$,
i.e., a hyperbolic structure on the complement of a finite set of distinguished points or punctures in $F$ together with the specification of a collection of horocycles in $F$, one about each puncture.  

Likewise,  (tangential) measures 
$\mu:\Delta\to{\mathbb R}$
give global linear coordinates on the space $\widetilde{\mathcal ML}_0(F)$ of decorated measured laminations in $F$ introduced and studied in
\cite{PapaPen}, where each puncture comes equipped with a foliated band of some width supplementing
a usual measured lamination in $F$ with compact support. In fact, $\mu= {\rm log}~ \lambda$, i.e.,
asymptotics in the sense of Thurston's boundary \cite{Thurston88} corresponds to taking logarithms.  
 This once-mysterious appearance of the logarithms is explained by tropicalization as follows.

Given any function 
$\phi(a,b,\cdots):{\mathbb R}_{>0}^n\to{\mathbb R}$, its {\it tropicalization} is defined as the limit as $t$ goes to infinity of ${1\over t}~{\rm log}~\phi(e^{ta},e^{tb},\cdots)$, if it exists,  taking values in the tropical semi-ring \cite{tropical}.  If a formula for $\phi$ is written by combining the coordinate functions of ${\mathbb R}^n$ using only multiplication and addition, its tropicalization is given by replacing each occurrence of multiplication by addition and each occurrence of addition by the binary maximum function with the precedence of tropical addition over maximum.  In particular, an inequality (either strict or weak) of two such expressions built only from multiplication and addition implies the corresponding (weak) inequality on their tropicalizations but not conversely.

We can nowadays say simply that \cite{PapaPen} showed that lambda lengths tropicalize to tangential measures,
a fact whose further ramifications are studied here.
Just as these results provided a key explicit example and gave clues
for general compactifications of cluster varieties \cite{FocGon},
so too we hope that the discoveries of this paper might inform
the tropicalization of extremal problems in related contexts,
viz. extending the cell decomposition of Riemann's moduli space
to general moduli spaces of flat connections.

Indeed, we study several more refined
notions regarding this procedure as follows.  There are two key conditions on  $\lambda:\Delta\to{\mathbb R}_{>0}$:

\vskip 2mm

\noindent \underline{Classical Triangle Inequality CT}: for every triangle $t\in T$ with frontier edges $x_i$, $i=1,2,3$, we have $$\lambda(x_i)<\lambda(x_j)+\lambda(x_k)~{\rm for}~\{ i,j,k\}=\{ 1,2,3\};$$

\vskip 2mm

\noindent \underline{Classical Face Condition CF}: In the notation of Figure \ref{fig:ptolemy} for the edges nearby
$e\in\Delta$, we have 
$$\lambda(e)\bigl(\lambda(a)\lambda(b)+\lambda(c)\lambda(d)\bigr)\leq \lambda(f)\bigl((\lambda(a)\lambda(d)+\lambda(b)\mu(c)\bigr),$$

\vskip 2mm

\noindent with respective tropical analogues for $\mu:\Delta\to{\mathbb R}$:

\vskip 2mm

\noindent  \underline{Tropical Triangle Inequality TT}: for every $t\in T$ with frontier edges $x_i$, $i=1,2,3$, we have $$\mu(x_i)\leq max \{ \mu(x_j),\mu(x_k)\}~{\rm for}~\{ i,j,k\}=\{ 1,2,3\};$$

\vskip 2mm

\noindent \underline{ Tropical Face Condition TF}: In the notation of Figure \ref{fig:ptolemy} for the edges nearby
$e\in\Delta$, we have $$\mu(e)+ max\{ \mu(a)+\mu(b),\mu(c)+\mu(d)\} \leq \mu(f)+ max \{ \mu(a)+\mu(d),\mu(b)+\mu(c)\}.$$

\vskip 2mm

Classically, CF  is equivalent to minimizing the objective function
$$\sum_{t\in T} \lambda(a)\lambda(b)\lambda(c),$$ where the sum is over all triangles $T$ complementary to $\Delta$, and where $t\in T$ has
frontier edges $a,b,c$.  This, in turn, is equivalent to a certain convexity condition \cite{Penner87,Penner11}
in Minkowski space, which is used to associate to $\Delta$ a corresponding cell $C(\Delta)\subset\widetilde{\mathcal T}(F)$
in the cell decomposition  of $\widetilde{\mathcal T}(F)$.  Furthermore, 
CF implies CT for positive lambda lengths.

Tropically, TF is equivalent (see Lemma \ref{L:corner}) to minimizing the objective
function given by the ordered list 
$${\rm sort}_\downarrow \{\mu(a)+\mu(b)+\mu(c):t\in T\}$$ of ``perimeters'' of triangles
as we shall see, and in fact, there is an algorithm (see Theorem \ref{T:seq})  to achieve the minimum contributed by Dylan Thurston.  Furthermore, if a global minimum also minimizes the total $\sum_{e\in\Delta} \mu(e)$, then it satisfies TT as well (see Corollary \ref{c:exist}).  This resolves a basic question about the existence of ideal triangulations appropriately suited to measured laminations and indeed provides an effective algorithm for finding them.

The condition TT is most interesting.  First of all, TT implies each of the other three conditions  (see Lemma \ref{implications}).  Secondly as a topological condition (see Lemma \ref{closedleaves}),
TT on tangential measure implies that every leaf is closed, i.e., the lamination is a multicurve possibly together with other puncture-parallel curves, i.e., it is a decorated weighted multicurve.

Another basic question is which collections of short curves are accessible from within a given cell in the decomposition of $\widetilde{\mathcal T}(F)$.
This was first answered in \cite{PenMc} relying on rather involved estimates on the asymptotics of lambda lengths
leading to a Fulton-MacPherson type treatment \cite{FM} of rates of their divergence.  The more general result on which
projectively measured laminations are accessible from within a given cell is answered here ``in a trice''
(see Theorem \ref{asymptotics-icd}) from the facts that
CF $\Rightarrow$ CT for positive lambda lengths and CT $\Rightarrow$ TT in the limit by tropicalization.
Just as in \cite{PenMc}, the proof is actually somewhat more involved in the non-generic case leading to the study of  balanced measured laminations (in Section \ref{cyclicsection}), which play the role for laminations that cyclic polygons play for geometric structures.

This paper is organized as follows.  In $\S$\ref{background}, we recall notation and facts about 
$\widetilde{\mathcal T}(F)$ and $\widetilde{\mathcal ML}_0(F)$ culminating in a theorem the author
proved with Athanase Papadopoulos in 1993
that lambda lengths on decorated Teichm\"uller space
tropicalize to measures
on decorated measured laminations
(though this result pre-dates tropical math by several years).
In $\S$\ref{S:trop}, we systematically study TT, first,
in relation to CT, CF and TF,
and second, a measured foliation satisfying TT with respect to some ideal triangulation is shown to 
have only closed leaves.  In $\S$\ref{cyclicsection}, we tropicalize two classical results \cite{Penner87,Penner11} on cyclic polygons leading to a new and naturally occurring class of measured laminations ``balanced'' with respect to an ideal cell decomposition, namely, for any completion to an ideal triangulation, TT holds on each triangle.  These tropical results are then applied in the following
section to give the measured lamination version of the earlier work with Greg McShane
on degenerations. 
An algorithm that plays the tropical role of the convex hull construction \cite{Penner87,Penner11}
is described in $\S$\ref{S:simp} and applied to show that for any weighted multicurve, there is an ideal triangulation so that the associated measure satisfies TT.

\section{Notations and Recollections}\label{background}

\vskip .2in

The {\it decorated Teichm\" uller space} $\widetilde {\mathcal T}=\widetilde {\mathcal T}(F)$ of a punctured surface
$F=F_g^s$ of genus $g$ with $s\geq 1$ punctures and negative Euler characteristic is the total space of
the trivial ${\mathbb R}_{>0}^s$-bundle $\pi : \widetilde {\mathcal T}\to {\mathcal T}$ over the usual Teichm\"uller space ${\mathcal T}={\mathcal T}(F)$ of $F$, where the fiber 
over a point of ${\mathcal T}$ is the set of all possible $s$-tuples of horocycles, one about each puncture  \cite{Penner87,Penner11}.  
Let $MC=MC(F)$ denote the mapping class group of $F$,  which acts naturally and $\pi$-equivariantly on ${\mathcal T}$ and $\widetilde {\mathcal T}$.  The total space $\widetilde {\mathcal T}$ admits a $MC$-invariant ideal cell decomposition as well as explicit global affine coordinates
which will be recalled presently.

The space 
$\mathcal{ ML}_0=\mathcal {ML}_0(F)$ 
of (possibly empty) measured laminations of compact support  on $F$ was introduced in \cite{Thurston88}.  The quotient $$\mathcal {PL}_0=\mathcal {PL}_0(F)=(\mathcal {ML}_0(F)-{\bf 0})/{\mathbb R}_{>0}$$ under the homothetic action of ${\mathbb R}_{>0}$ on transverse measures  is a PL-sphere of dimension $6g-7+2s$ called {\it Thurston's boundary}, where ${\bf 0}$ denotes the empty lamination, naturally compactifying ${\mathcal T}\approx {\mathbb R}^{6g-6+2s}$ to a closed ball 
$$\overline{{\mathcal T}}=\overline{{\mathcal T}}(F)={\mathcal T}(F)\sqcup {\mathcal{PL}}_0(F),$$
 upon which $MC$ acts naturally and continuously.
Let $ [{\mathcal L}]\in \mathcal {PL}_0$ denote the projective 
class of $\mathcal F\in\mathcal {ML}_0-\{ {\bf 0}\}$.

Following \cite{PapaPen}, there is furthermore a vector space
$\widetilde{\mathcal{ML}}_0=\widetilde{\mathcal{ML}}_0(F)$, the space {\it decorated measured laminations}, together with a continuous map $\pi:\widetilde{\mathcal{ML}}_0\to{\mathcal {ML}}_0$, where the fiber over a point of $\mathcal{ML}_0$ is isomorphic to ${\mathbb R}^s$.  A decorated measured lamination is regarded as a measured lamination with compact support in the usual sense except that puncture-parallel leaves, called {\it collar curves}, are now allowed and are given a 
${\mathbb R}$-valued ``decoration''; if the decoration is positive, then this parameter is taken as the transverse measure of a foliated band parallel to the collar curve, and if it is negative, then it is thought of as a ``deficit'' of such leaves as in \cite{PapaPen}.
In particular, the fiber over $\bf 0$ is the collection of all ${\mathbb R}$-weighted collar curves.

 An {\it ideal triangulation} $\Delta$ of $F$ is (the isotopy class of) a collection of arcs decomposing $F$ into triangles with their vertices at the punctures.  A {\it flip} on an arc $a$ in $\Delta$ which separates distinct triangles complementary to $\Delta$ removes $a$ from $\Delta$ and replaces it with the unique distinct edge.  The groupoid $Pt=Pt(F)$ whose objects are ideal triangulations and whose morphisms are flips between them is called the {\it Ptolemy groupoid}.  The groupoid
$M\Gamma=M\Gamma(F)$ whose objects are $MC$-orbits of ideal triangulations imbued with an enumeration of the edges, and whose morphisms are $MC$-orbits of pairs of such with the obvious mutliplication is called
the {\it mapping class groupoid}.  $Pt$ and $M\Gamma$ admit simple presentations \cite{Penner87,Penner11}, and $M\Gamma$ contains $MC$ as the stabilizer of any object, which is a subgroup of finite index.

\begin{figure}[!h]
\begin{center}
\epsffile{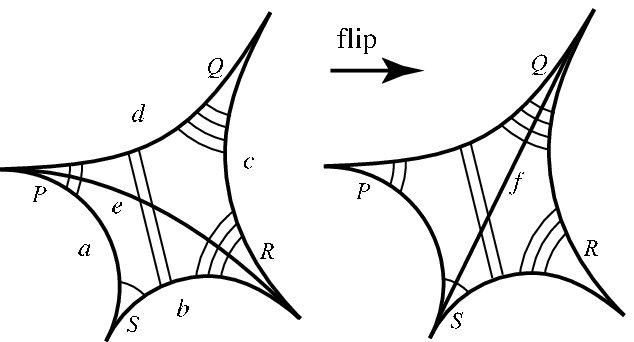}
\caption{A flip}
\label{fig:ptolemy}
\end{center}
\end{figure}

Upon fixing an ideal triangulation $\Delta$, there are global affine coordinates $c_L :{\mathbb R}_{>0}^\Delta\to \widetilde T$ on decorated Teichm\"uller space called {\it lambda lengths} or ``Penner coordinates'' \cite{Penner87,Penner11}.  The effect of a flip
on lambda lengths is described by a {\it Ptolemy transformation} $ef=ac+bd$ in the notation of 
Figure~\ref{fig:ptolemy}, where we identify an arc with its lambda length for convenience; furthermore, the cross-ratio of the quadrilateral illustrated is given by $ac/bd$ in lambda lengths.

There are moreover global linear coordinates
$c_M:{\mathbb R}^\Delta\to \widetilde{\mathcal{ML}}_0$ called {\it (tangential) measures} on $\Delta$ regarded as a train track with stops, cf. \cite{PapaPen,PenHar}.  The effect of a flip on measures is described by $e+f=max \{ a+c,b+d\}$ in the notation of Figure \ref{fig:ptolemy}, where now we identify an arc with its measure for convenience.

\begin{figure}[!h]
\begin{center}
\epsffile{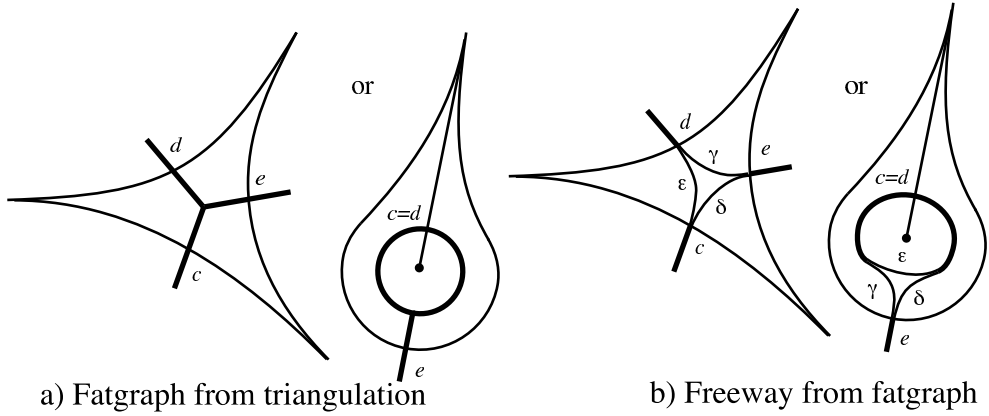}
\caption{Fatgraphs and freeways}
\label{fig:freeway}
\end{center}
\end{figure}

For each ideal triangulation $\Delta$, there is an associated branched one-submanifold called the {\it dual freeway} $\tau =\tau_\Delta$ as illustrated in Figure \ref{fig:freeway}b, where there is a small triangle of $\tau$
lying inside each complementary region to $\Delta$ and one dual ``long'' edge meeting each arc of $\Delta$, cf. \cite{PapaPen}.
The measure on the train track with stops $\Delta$ can dually be regarded as a transverse measure on the freeway $\tau$; the values on the long edges of $\tau$ are unrestricted real numbers which extend uniquely to a ${\mathbb R}$-valued transverse measure on all of $\tau$, which in turn uniquely determines a decorated measured lamination. 
Notice that in the notation of Figure~\ref{fig:freeway} on a freeway inside a triangle, we have $\gamma= {{d+e-c}\over 2}$,
$\delta ={{c+e-d}\over 2}$ $\varepsilon={{c+d-e}\over 2}$, and inside a monogon, we have
$\gamma=\delta={e\over 2}$, $\varepsilon=c=d$.

An {\it ideal cell decomposition} of $F$ is a subset of an ideal triangulation each of whose complementary regions is simply connected.
The dual of an ideal cell decomposition is a {\it fatgraph} spine embedded in $F$ as illustrated for the generic case of ideal triangulations in Figure
\ref{fig:freeway}a.

There is a $MC$-invariant ideal cell decomposition  \cite{Harer85,Penner87,Penner11,Strebel} of $\widetilde {\mathcal T}$ whose cells are in one-to-one correspondence with ideal cell decompositions of $F$, where the face relation is induced by inclusion.  
The open cell $C(\Delta)\subseteq \widetilde {\mathcal T}$ associated to an ideal triangulation $\Delta$ is determined  by the positivity $E>0$ of a collection of expressions
$$E={{a^2+b^2-e^2}\over{abe}}~+~{{c^2+d^2-e^2}\over{cde}}$$
in the lambda length notation of Figure~\ref{fig:ptolemy} called {\it simplicial coordinates}, one for each edge of $\Delta$.    These
simplicial coordinates play a key role in this paper.

Moreover,
the cell $C(\Delta')$ corresponding to an ideal cell decomposition $\Delta'$ lies in the closure of $C(\Delta)$ if and only if $\Delta'\subseteq\Delta$; in this case, $C(\Delta')$ is determined by the vanishing
of simplicial coordinates of edges of $\Delta-\Delta'$, and non-vanishing putative simplicial coordinates on the edges of an ideal cell decomposition $\Delta'$ uniquely determine lambda lengths on the edges of any ideal triangulation containing it (as will be further explained in Section \ref{cyclicsection}). 

In particular, the simplicial coordinates are the hyperbolic analogues of Strebel coordinates in the conformal setting, but there is no conformal analogue of lambda lengths.  
The explicit calculation of lambda lengths from simplicial coordinates is the ``arithmetic problem'' discussed at length in \cite{Penner95,Penner11}.

There is the following diagram summarizing the maps and spaces just defined:

\vskip .2cm

$$\begin{array}{cccccccccccccccc}
&&&&c_L&&&&~~~{\rm log}&&&&c_M\\
&&\widetilde {\mathcal T}&&\leftarrow&&{\mathbb R}_+^\Delta&&~~~{ \rightarrow}&&{\mathbb R}^\Delta&&\rightarrow&&\widetilde{\mathcal{ML}}_0\\
&&&&&&&& \\
&&~~~\downarrow~\pi&&&&&&&&&&&&\downarrow~\pi\\\\
&&{\mathcal T}&&&&&&&&{\mathcal{ML}}_0-\{{\bf 0}\}&&\subset&&{\mathcal{ML}}_0\\\\
&&\cap&&&&&&&&\downarrow~[\cdot]\\\\
&&\overline {\mathcal T}&&&&\supset&&&&\mathcal{PL}_0\\
\end{array}$$

\vskip .2cm

\noindent where ${\rm log}:{\mathbb R}_+^\Delta\to{\mathbb R}^\Delta$ denotes the function which in each coordinate simply takes natural logarithms, $c_L$ denotes lambda lengths and $c_M$ denotes measures on $\Delta$, or equivalently transverse measures on $\tau$.

~~\vskip .2cm

\begin{theorem}\label{theoremone}\cite{PapaPen,Penner92,Penner11}
~Fix an ideal triangulation $\Delta$ of $F$ and consider a sequence $\lambda_i\in{\mathbb R}_+^\Delta$, for $i\geq 0$, escaping each compactum in Teichm\"uller space.

\vskip .2cm

\noindent Then $[\pi\circ c_L(\lambda_i)]\to [{\mathcal L}]\in {\mathcal {PL}}_0$ in the topology of Thurston's boundary $\overline{{\mathcal T}}$ of Teichm\"uller space ${\mathcal T}$ if and only if $\bar \lambda_i\to\mu\in{\mathbb R}^\Delta$, where
$$\pi\circ c_M(\mu )={\mathcal L}\in\mathcal{ML}_0.$$

\noindent Furthermore, the tropicalization $$\bar e+\bar f  = max \{ \bar a+\bar c,\bar b+\bar d\}$$ of the Ptolemy equation $ef=ac+bd$ describes the effect of flips on measures,
where we write $\bar x={\rm log}~x$ for $x\in{\mathbb R}_+$.

\end{theorem}

\begin{remark}
Moreover,  the expression
$$2\sum d\bar a\wedge d\bar b +d\bar b\wedge d\bar c + d\bar c\wedge d\bar a,$$
where the sum is over all triangles complementary to $\Delta$ with coordinates $a,b,c$ in this clockwise cyclic order corresponding to the boundary describes the pull-back 
\cite{Penner92,Penner11} by $\pi\circ c_L$ of the Weil-Petersson K\"ahler two-form on ${\mathcal T}$  in logarithms $\bar a, \bar b, \bar c$ of lambda lengths as well as the pull-back \cite{PapaPen}
by $\pi\circ c_M$ of the Thurston symplectic form on $\mathcal{ML}_0$ in measures $\bar a,\bar b,\bar c$.
\end{remark}

\vskip .3cm

\section{Tropical Triangle Inequalities}\label{S:trop}

\vskip .2cm
 
Fix an ideal triangulation $\Delta$ of $F$ and an assignment of real numbers called ``weights'', one to each arc in $\Delta$, and consider the following inequalities on weights, where we shall often identify an edge with its weight for convenience.

\vskip .3cm

\noindent {\bf (Classical) Triangle Inequalities CT}: for each triangle complementary to $\Delta$ whose boundary has consecutive edges $a,b,c$, we have $$\aligned
a&< b+c,\\
b&< a+c,\\
c&< a+b;\\
\endaligned$$

\vskip .2cm

\noindent{\bf Tropical Triangle Inequalities TT}: in the notation of the triangle inequalities, the maximum $max \{ a,b,c\}$ is achieved either twice or thrice, that is, 
$$\aligned a&\leq  max \{ b,c\},\\
b&\leq max \{ a,c\},\\
c&\leq max \{a,b\};\\
\endaligned$$

\vskip .2cm

\noindent {\bf (Classical) Face Conditions CF}: for each quadrilateral comprised of two adjacent triangles complementary to $\Delta$ whose edges are as illustrated in Figure~\ref{fig:ptolemy}, we have
$$\bigl [ cd(a^2+b^2-e^2)\bigr ] + \bigl [ ab(c^2+d^2-e^2)\bigr ]\geq 0,$$ i.e., we have $e^2(ab+cd)\leq (ac+bd)(ad+bc)$, or in other words by the Ptolemy equation, we have $$e(ab+cd)\leq f(ad+bc);$$

\vskip .2cm

\noindent{\bf Tropical Face Conditions { TF}}: in the notation of the face conditions, we have
$$\aligned
2e+max \{ a+b,c+d\} &\leq max\bigl\{ c+d+2 max \{ a,b\},a+b+2 max \{ c,d\}\bigr \}\\
&=max \{ a+c,b+d
\}+max \{ a+d,b+c\},\\
\endaligned$$
where the equality follows as the
tropicalization of the classical identity $ab(c^2+d^2)=(ac+bd)(ad+bc)$;
in other words by the tropical Ptolemy equation, we have
 $$e+max \{ a+b,c+d\} \leq f+max \{ a+d,b+c\}.$$
 
 \begin{lemma}
 Fix an ideal triangulation $\Delta$ of the surface $F$.  A collection $\lambda\in{\mathbb R}_{>0}^\Delta$
 of lambda lengths describes a point in the corresponding open cell $C(\Delta)\subseteq{\mathcal T}(F)$
 if and only if $\lambda$ minimizes the objective function $\sum_{t\in T} \lambda(a)\lambda (b)\lambda (c)$, where the sum is over all triangles $t$ complementary to
 $\Delta$ with frontier edges $a,b,c$, i.e., no flip on edges of $\Delta$ decreases this objective function.
 \end{lemma}
 
 \begin{proof}
 This follows directly from the last formulation of CF.
 \end{proof}

\vskip .3cm

 In effect in Section \ref{S:simp}, we shall learn that the natural tropicalization
of this objective function is the ordered sequence of perimeters
$${\rm sort}_\downarrow\{\mu(a)+\mu(b)+\mu(c):t\in T~{\rm with~frontier~edges}~a,b,c\},$$ for $\mu\in{\mathbb R}^\Delta$, 
in the sense that CF is equivalent to minimizing the former objective function for decorated hyperbolic structures and TF is equivalent to minimizing the latter for decorated measured laminations. Insofar as CF is equivalent to convexity in a convex hull construction in Minkowski space \cite{Penner87,Penner11}, so too we think of TF as a kind of tropical convexity condition.

There are various implications among these conditions:

\begin{lemma}\label{implications} Fix a weight on an ideal triangulation.

\vskip 2mm

\noindent {\rm TT} implies {\rm CT, CF}, and {\rm TF}, and these are the only implications that hold among these conditions.

\vskip 2mm

\noindent  {\rm CF} implies {\rm CT} for positive weights.

\vskip 2mm

\noindent  {\rm TF} implies {\rm TT} provided the sum of the all the weights is minimal among such triangulations.

\end{lemma}

\begin{proof}
The non-uniqueness of a maximum weight on edges of triangles complementary to an ideal triangulation $\Delta$ for tropical triangle inequalities shows that the classical triangle inequalities hold not only on the weights themselves, but also their squares.  It follows that TT $\Rightarrow$ CT and also that TT $\Rightarrow$ CF. To see that TT $\Rightarrow$ TF in the notation of Figure \ref{fig:ptolemy}, we have
$e\leq max \{ a,b\}$ and $e\leq max \{ c,d\}$ in particular, so
$$\aligned
2e+max \{ a+b,c+d\}&= max \{ a+b+2e,c+d+2e\} \\
&\leq max \bigl \{ a+b+2max \{ c,d\}, \\
&\hskip 1.6cm c+d+2max \{ a,b\}\bigr\}\\
&=max \{ a+c,b+d\} +max \{ a+d,b+c\},\\
\endaligned$$
as claimed, where the last equality follows as before.

To see that these are the only implications that hold among the conditions on weights, it remains to give a series of counter-examples $a,b,c,d,e$ as usual in the notation of Figure \ref{fig:ptolemy}  as follows:

\vskip .2cm

\noindent 4,4,1,6,5 shows TF $\not\Rightarrow$ TT; 4,3,3,4,5 shows CT $\not\Rightarrow$ TT; 
3,2,3,3,4 shows CT $\not\Rightarrow$ CF; 3,1,1,1,1 shows TF $\not\Rightarrow$ CT;
and 2,2,4,4,3 shows CF $\not\Rightarrow$ TF.

\vskip .2cm

\noindent Extending these specifications to weights defined on an entire ideal triangulation is an exercise, and likewise is the verification that these counter-examples together with the positive results determine all implications among these conditions.

\vskip .2cm

For the second assertion that CF $\Rightarrow$ CT for positive weights, we argue as in \cite{Penner87,Penner11} as follows. Adopt the notation of Figure \ref{fig:ptolemy} in the universal
cover if the edges $a,b,c,d$ are not distinct or if $e$ does not have distinct endpoints.  If
$c+d\leq e$, then $c^2+d^2-e^2\leq -2cd$.  The face condition 
$0\leq cd(a^2+b^2-e^2)+ab(c^2+d^2-e^2)$ then
gives $0\leq cd  [(a-b)^2-e^2]$, and we find a second triangle so that the triangle inequality fails.
It follows that any failure of the triangle inequalities implies that there is a cycle of triangles of such
failures, namely, there are complementary triangles $t_0,\ldots ,t_n=t_0$, where the frontier of 
$t_i$ consists of arcs $e_{i-1},e_i,b_i$ with $e_i=t_i\cap t_{i+1}$, for $i=0,\ldots ,n$ and the subscripts are taken modulo $n$, so that $e_j\geq b_j+e_{j-1}$, for $j=1,\ldots, n$.
Upon summing and canceling
like terms, we find $0\geq \sum _{j=1}^n b_j$, which is 
absurd for positive weights.

For the proof of the final assertion that TF implies TT under suitable conditions, suppose to the contrary that there is a triangle not satisfying TT with an edge $e$.  Take the maximal weight such edge $e$ and adopt the nearby notation as in Figure \ref{fig:ptolemy}.  Let $p,q,r,s$ denote the widths of the bands parallel to
$P,Q,R,S$
of a decorated measured lamination $\mu\in{\mathbb R}^\Delta$ 

TF evidently implies either $p$ or $r$ is minimum among $p,q,r,s$ (cf. also Lemma \ref{L:corner}), and we may suppose without loss that it is $p$. If $\mu(b)>\mu (e$), then we must have $\mu(a)=\mu(b)$, which implies that  there are no leaves of $\mu$ in the quadrilateral near $e$ other than those parallel to the cusps and furthermore $p=r$, which in turn implies that $\mu(c)=\mu(d)$ as well.  Similar remarks apply in case $\mu(c)$ or $\mu(d)$ is greater than $\mu (e)$  Thus, $\mu(e)$ must be maximal among its neighbors.

Now, since $\mu(e)$ is maximal, then in particular $\mu(e)\geq \mu(b)$, i.e., letting $x$ denote width of the band of leaves
that cross the two sides of quadrilateral near $e$, we have $p+x+r\geq s+x+r$, that is, $p\geq s$.  Since we also necessarily have $p\leq s$, we conclude $p=s$.  Thus, the flip here is neutral with respect to TF.

Furthermore, $\mu(e)\geq \mu (d)$ gives $p+x+r\geq p+x+q$, that is $r\geq q$.  If this is a strict inequality,
then $\mu(f)=s+x+q=p+x+q<p+x+r=\mu (e)$, so we can flip the edge $e$ and decrease the total weight, which is contrary
to hypothesis.
\end{proof}

\vskip .2cm

Consider a positive measure on the frontier edges of a triangle $t$ satisfying TT.  The corresponding measured lamination of $t$ consists of three bands of leaves of common width $a$ connecting each pair of frontier edges of $t$ together with exactly one further band of width $b$ connecting a single pair of frontier edges of $t$, for some $a,b\geq 0$.  Of course, a positive measure on an ideal triangulation $\Delta$ of a surface $F$ must restrict to such a special measured lamination on each triangle complementary to $\Delta$, and TT imposes further global constraints as follows.

\vskip .2cm

\begin{lemma}\label{closedleaves}
Suppose that $\mu:\Delta\to{\mathbb R}$ is the measure of a decorated measured lamination 
with respect to some
ideal triangulation $\Delta$ and that 
$\mu$ satisfies the tropical triangle inequalities on each triangle complementary to $\Delta$.
Then each leaf of the corresponding lamination ${\mathcal L}=c(\mu)$ is a simple closed curve.
\end{lemma}

\begin{proof}
The proof is by downward induction on the number of times that $\mu$ achieves its maximum value. If this number agrees with the number of edges of $\Delta$, i.e., if all measures are equal, then the corresponding decorated measured lamination is a collection of equally weighted collar curves by construction. 

For the induction, suppose that the lemma is established if the measure takes its maximum value $n+1$ times, and consider a measure $\mu$ taking its maximum value $n$ times.
Let $m$ be the maximum value of $\mu$ and say that an arc in the ideal triangulation $\Delta$ is  \textit{maximal} if its measure is exactly $m$. 
Split each band of leaves of ${\mathcal L}=c(\mu)$ crossing a maximal edge into three consecutive bands of respective widths ${1\over 2} (m-\varepsilon)$, $ \varepsilon$, and ${1\over 2}(m-\varepsilon)$.  

We claim that for sufficiently small $\varepsilon$, if a leaf of ${\mathcal L}$ belongs to the central band for one edge, then it will intersect maximal edges only and belong to the central band for each of them. 
To see this, follow a leaf from a central band.  It obviously stays in the central band when crossing a triangle with three maximal edges for any value of $\varepsilon$. 
When the leaf crosses a triangle with two maximal edges and one edge of measure $m'< m$, it evidently stays in the central band if and only if $\varepsilon \leq m-m'$.  In particular, the leaves in central bands are thus simple closed curves for sufficiently small $\varepsilon$.

Let us take the maximum possible $\varepsilon$ satisfying this condition, namely, the minimum over all triangles with exactly two maximal edges of the quantity $m-m'$, where $m'<m$ is the measure of the edge of the triangle which is not maximal.
If we remove from ${\mathcal L}$ the sub-lamination consisting of all of the simple closed leaves arising from central bands to produce a measured lamination ${\mathcal L}'$, then the measure of 
${\mathcal L}'$ on $\Delta$ again satisfies the tropical triangle inequalities with the maximal edges of measure $m-\varepsilon$. 

On the other hand, since at least one non-maximal edge has exactly this measure by our choice of $\varepsilon$, the maximal value is achieved more than $n$ times for ${\mathcal L}'$ .  All leaves of ${\mathcal L}'$ are therefore simple closed curves by the induction hypothesis, and so  ${\mathcal L}$ has this property as well.
\end{proof}

\begin{corollary}\label{three}
Suppose that $\mu:\Delta\to{\mathbb R}$ is an integral measure 
satisfying the tropical triangle inequalities on an ideal triangulation $\Delta$. 
If $\mu$ corresponds to a simple closed curve, then its value on each edge of $\Delta$ is at most three. It follows that
for each fixed ideal triangulation $\Delta$, there are only finitely many simple closed curves whose corresponding measures on $\Delta$ satisfy the tropical triangle inequalities.
 \end{corollary}

\begin{proof}
The first part follows from the proof of the previous result where the maximum possible $\varepsilon$ is at least one for an integral measure, and the second part follows immediately from the first.
\end{proof}

\vskip .3cm

\section{Cyclic polygons and balanced laminations}\label{cyclicsection}

\vskip .2cm

A tuple $a_1,\ldots ,a_n$ of real numbers might satisfy the following conditions:

\vskip .2cm

\noindent {\bf (Classical) Generalized Triangle Inequalities CGT}:~for each $i=1,\ldots, n$, we have $a_i< \sum_{j\neq i} a_j$;

\vskip .2cm

\noindent{\bf Tropical Generalized Triangle Inequalities TGT}:~the maximum among $\{ a_i\}_1^n$
is achieved at least twice, i.e., for each $i=1,\ldots ,n$, we have
$a_i\leq max \{ a_j:j\neq i\}$.

\vskip .2cm

\begin{lemma}\label{CTCGT}
Fix an ideal triangulation $\Delta$ of a polygon $P$ and consider a collection of weights on the arcs in $\Delta$ including the frontier arcs of $P$.  Then CT on the weights on $\Delta$ implies CGT on the
weights on the frontier of $P$, and TT likewise implies TGT on the frontier.
\end{lemma}

\begin{proof}
Each assertion follows by induction on the number of sides of $P$ with the vacuous basis
step for $P$ a triangle.
\end{proof}

\vskip .2cm

We shall say that a Euclidean polygon is {\it cyclic} if it inscribes in a circle and recall the following two results:

\vskip .2cm

\begin{lemma}\label{cyclicquad}
Consider a cyclic Euclidean polygon with a fixed ideal triangulation $\Delta$.  If $e\in\Delta$ with nearby edges $a,b,c,d$ as depicted in Figure \ref{fig:ptolemy}, then the Euclidean lengths are related by $$e^2=
{{(ac+bd)(ad+bc)}\over{(ab+cd)}}.$$ 
\end{lemma}

\begin{proof}
This is precisely the case of vanishing of simplicial coordinate in the face condition.
\end{proof}

\begin{theorem} \label{cyclic} \cite{Penner87,Penner11} A tuple $a_i\in{\mathbb R}_{>0}$, for $i=1,\ldots ,n$, occurs as the edge lengths
of a cyclic Euclidean polygon if and only if it satisfies the classical generalized triangle inequalities {\rm CGT}.  Furthermore,
such tuples of lengths of frontier edges of polygons uniquely determine the corresponding cyclic Euclidean polygon.
\end{theorem}

\vskip .2cm

In this section, we give the tropical analogues of Lemma \ref{cyclicquad} in Lemma \ref{balancedlemma}
and of Theorem \ref{cyclic} in Proposition \ref{balanced}.

If $\Delta'$ is an ideal cell decomposition and $\Delta\supseteq\Delta'$ is an ideal triangulation of some surface $F$, then a point of $C(\Delta')\subseteq \widetilde {\mathcal T}$ is completely determined by the lambda lengths of the arcs in $\Delta'$ as follows.  Fix some complementary region $R$ of $\Delta'$.  The lambda lengths on the frontier of $R$ satisfy CGT by Lemmas \ref{implications} and \ref{CTCGT} and so determine a cyclic Euclidean polygon according to Theorem \ref{cyclic}.  An edge of $\Delta - \Delta'$ triangulating $R$ has its lambda length given by the Euclidean length of the corresponding diagonal of this cyclic polygon.  In particular, lambda lengths on the frontier edges of a cyclic polygon properly embedded in a decorated surface uniquely determine the lambda lengths on all diagonals of this polygon.

We shall say that a measured lamination of a polygon is {\it balanced} provided that for each ideal triangulation of it, the tropical triangle inequalities hold on each complementary triangle.  

\begin{lemma}\label{balancedlemma}
Consider a polygon with fixed ideal triangulation $\Delta$.  If $e\in\Delta$ with nearby edges $a,b,c,d$ as depicted in Figure \ref{fig:ptolemy}, then the measures of any balanced measured lamination of the polygon are related by $$e={1\over 2}\biggl \{
max \{ a+c,b+d\} + max \{ a+d,b+c\} -max \{ a+b,c+d\}\biggr \}.$$ 
Furthermore, $a,b,c,d\geq 0$ implies $e\geq 0$.
\end{lemma}

\begin{proof}
According to Lemma \ref{implications}, the assumed tropical triangle inequalities
imply the tropical face conditions on $\Delta$, and likewise for the ideal triangulation that arises from 
$\Delta$ by a flip on $e$.
It therefore follows that 
$$
\aligned
e&\leq {1\over 2}\bigl \{max \{ a+c,b+d\} +\max \{ a+d,b+c\}-max \{ a+b,c+d\} \bigr\},\\
f&\leq {1\over 2}\bigl \{max \{ a+c,b+d\} +\max \{ a+b,c+d\}-max \{ a+d,b+c\} \bigr\},\\
\endaligned$$
so the sum provides the estimate $e+f\leq max \{ a+c,b+d\}$.  On the other hand, the tropicalization of the
Ptolemy equation gives the equality $e+f=max \{ a+c,b+d\}$, and therefore equality must hold also in the estimates
above on $e$ and $f$.

To prove the final assertion about non-negativity, observe that $a,b,c,d$ satisfies TGT by Lemma \ref{CTCGT} since $a,b,e$ and $c,d,e$ satisfy
TT by hypothesis.  There are several cases depending upon which of $a,b,c,d$ achieve the 
maximum $max \{ a,b,c,d\}$:

\vskip .2cm

\noindent if $a=b$ is the maximum, then $e=max \{ c,d\}$; if $a=c$ or $a=d$ is the maximum, then
$e=a$; if $b=c$ or $b=d$ is the maximum, then $e=b$; and if $c=d$ is the maximum, then $e=max \{ a,b\}$.
\end{proof}

\begin{proposition}\label{balanced}
A tuple $a_i\in{\mathbb R}_{>0}$, for $i=1,\ldots ,n$, occurs as the measure on the frontier edges of a polygon $P$ for a balanced measured lamination if and only if it satisfies the TGT inequalities. 
Furthermore, such tuples of measures on the frontier edges of $P$ uniquely determine the corresponding measured lamination on $P$.
\end{proposition}

\begin{proof}
Necessity of TGT follows from Lemma \ref{CTCGT}.
For sufficiency and to complete the proof, we must show that there is a unique balanced measured lamination of $P$ realizing any given tuple of putative measures on its frontier edges which satisfies TGT.  The proof proceeds by induction on the number $n$ of sides of $P$.  The case  $n=3$, where balanced $\Leftrightarrow$ TT $\Leftrightarrow$ TGT, was discussed already in Section \ref{background}, and the case $n=4$, follows from the previous lemma.   

For the inductive step with $n\geq 5$, consider two consecutive frontier edges $a,a'$ of $P$ and the diagonal $d$ of $P$ separating them from the rest.  If the putative measures $a$ and $a'$ are different, then $d=max \{ a,a'\}$ is uniquely determined from the balanced condition.  Furthermore, suppose that $a,a',a''$ are consecutive frontier edges of $P$ with $d$ and $d'$, respectively, the diagonals separating $a,a'$ from $a''$ and separating $a$ from $a',a''$, and that $d''$ is the diagonal of $P$ separating $a,a',a''$ from the rest of the frontier.  If $a=a'=a''$, then
$d+d'=max \{ a+a'', a'+d''\}=a+max \{a,d''\}$ by the tropical Ptolemy equation, and so
$d<a$ implies that $d'>max \{ a,d''\}\geq a$, which violates the balanced condition.  
It therefore follows that $d=d'=a$ in this case.  In either case, cutting on the diagonal $d$ produces an $(n-1)$-gon to which we may apply the inductive hypothesis and thereby complete the proof of existence. 

As to uniqueness, any diagonal $d$ of $P$ separates the frontier edges of $P$ into two disjoint sets $A,A'$.  The determination of measure on $d$ in the previous paragraph is given by
$d=min\{ max A, max A'\}$ independently of any choices as desired.
\end{proof}

\vskip .2cm

\begin{remark}
In particular, measures of a balanced decorated measured lamination on the frontier edges of a polygon properly embedded in a surface uniquely determine the measures on all diagonals of this polygon.  According to Theorem \ref{theoremone}, 
these thus computable measures describe the more elusive asymptotics of the logarithms of lambda lengths.
\end{remark}

\vskip .3cm

\section{Asymptotics of cells}\label{S:asym}

\vskip .2cm

In this section, we give the tropical analogue of  \cite{PenMc}, namely, we answer the question of which measured laminations are the asymptotes of paths
of geometric structures within a fixed cell $C(\Delta)$ of the cell decomposition of $\tilde {\mathcal T}$.

\begin{lemma}\label{asymptotics-it}
Fix an ideal triangulation $\Delta$ and consider a sequence $\lambda_i\in{\mathbb R}_{>0}^\Delta$ with
$c(\lambda _i)\in C(\Delta)\subseteq \widetilde {\mathcal T}$.  Then $\pi\circ c(\lambda_i)\in {\mathcal T}$ is asymptotic to
$[\pi\circ c(\mu)]\in{\mathcal PL}_0$ in the topology of Thurston's boundary, for some $\mu\in{\mathbb R}^\Delta$, if and only if
$\mu$ satisfies the tropical triangle inequalities on $\Delta$.  In particular,  if $\pi\circ c(\lambda_i)$ tends to 
$[\pi\circ c(\mu)]$, where $c(\lambda_i)\in C(\Delta)$, then each leaf of $\pi\circ c(\mu)$ is a simple closed curve.
\end{lemma}

\begin{proof}
To prove the first assertion, suppose to begin
that  $\pi\circ c(\lambda_i)$ is asymptotic to $[\pi\circ c(\mu)]$.  According to Theorem \ref{theoremone}, the projective class of
the measure $\mu$ is given by the limiting values of logarithms of lambda lengths ${\rm log}~\lambda_i$ on each edge of $\Delta$.  Furthermore, each $\lambda_i>0$ satisfies the face conditions since $c(\lambda _i)\in C(\Delta)$ by hypothesis and hence also satisfies the triangle inequality by
Lemma \ref{implications}.  Since an inequality between expressions constructed using only addition and multiplication 
implies its tropicalization, it follows that the tropical triangle inequalities on measures must hold as well.  
The last assertion then follows directly from Lemma \ref{closedleaves}.  

For the converse, given $\mu\in{\mathbb R}^\Delta$ satisfying the tropical triangle inequalities, define
the one-parameter family of lambda lengths 
$$\lambda _t(e)=t^{\mu (e)},~~{\rm for}~e\in\Delta~{\rm and}~t\in{\mathbb R}_{>0}.$$
Consider an edge $e\in\Delta$ with nearby edges as illustrated in Figure~\ref{fig:ptolemy}.  
As before, since the maximum of $\mu$ is achieved either twice or thrice on $\{ a,b,e\}$ and on $\{ c,d,e\}$, each
summand in brackets in the face condition for $\lambda_t$ is strictly positive, and hence so too is their sum.  It follows that
$\lambda _t\in C(\Delta)$, for each $t\geq 1$, and according to Theorem \ref{theoremone}, we furthermore
have $\pi\circ c(\lambda _t)\in T$ is asymptotic to $[\pi\circ c(\mu)]\in{\mathcal PL}_0$ as $t$ goes to infinity as required.
\end{proof}

\vskip .2cm

\begin{theorem}\label{asymptotics-icd}
Fix an ideal triangulation $\Delta$ 
and consider a sequence $\lambda_i\in{\mathbb R}_{>0}^\Delta$ with
$c(\lambda _i)\in C(\Delta)\subseteq \widetilde {\mathcal T}$.  
Then $\pi\circ c(\lambda_i)\in {\mathcal T}$ is asymptotic to
$[{\mathcal L}]\in{\mathcal PL}_0$ in the topology of Thurston's boundary
if and only the restriction of ${\mathcal L}$ to each component of $F-\Delta'$
is balanced.  In particular,  if $\pi\circ c(\lambda_i)$ tends to 
$[{\mathcal L}]$, where $c(\lambda_i)\in C(\Delta)$, then each leaf of ${\mathcal L}$ is a simple closed curve.  Only finitely many simple closed curves, namely, those
whose measures on $\Delta$ are bounded above by three and satisfy the tropical triangle inequalities, can occur as leaves
\end{theorem}

\begin{proof}The generic case of an ideal triangulation $\Delta '=\Delta$ is covered by the previous result, which furthermore shows that the balanced condition is necessary since $C(\Delta ')\subseteq C(\Delta)$.

As to sufficiency, let $\mu:\Delta \to{\mathbb R}_{>0}$ denote the measure of a measured lamination ${\mathcal L}$ which is balanced on each component of $F-\Delta'$.  Define a one-parameter family of lambda lengths 
$$\lambda _t(e)=t^{\mu (e)},~~{\rm for}~e\in\Delta'~{\rm and}~t\in{\mathbb R}_{>0}.$$
Again by Lemmas \ref{implications} and \ref{CTCGT}, $\lambda _t$ satisfies CGT on components of
$F-\Delta'$, for each $t$,  and so uniquely determines a cyclic Euclidean polygon by Theorem \ref{cyclic} corresponding to each such component.  If $e\in\Delta-\Delta'$, then let $\lambda _t(e)$ be given by the Euclidean length of the corresponding diagonal of the associated cyclic polygon.  This defines
$\lambda _t$ on each edge of $\Delta$ for each $t$.

By construction, we have $\pi\circ c(\lambda _t)\in C(\Delta')$, for all $t$, and by Lemma \ref{asymptotics-it}, $\pi\circ c(\lambda _t)\in C(\Delta')$ is asymptotic to a measured lamination ${\mathcal L}'$ whose measure $\mu '$ on $\Delta$ agrees with $\mu$ on $\Delta'\subseteq\Delta$.  Since the lambda lengths $\lambda_t$ satisfy CT on $\Delta$, $\mu '$ must satisfy TT on $\Delta$.  Similarly, for any other ideal triangulation $\Delta_1\subseteq\Delta'$, the analogously defined lambda lengths on $\Delta_1$ describe the same  path in  $\widetilde {\mathcal T}$, therefore have the same asymptote ${\mathcal L}'$, and likewise have measure on $\Delta_1$ satisfying TT.
It follows that ${\mathcal L}'$ is balanced on each component of $F-\Delta'$.

Thus, ${\mathcal L}$ and ${\mathcal L}'$ are both balanced on each component of $F-\Delta'$.  Since their measures $\mu$ and $\mu'$ agree on $\Delta'$ by construction, we must have ${\mathcal L}={\mathcal L}'$  by the last part of Proposition \ref{balanced}.

The last part follows from Corollary \ref{three}.
\end{proof}

\vskip .3cm

\section{Simplifying Triangulations}\label{S:simp}

Fix any ideal triangulation $\Delta$ of a punctured surface $F$ with complementary triangles $T$,
and consider a (decorated) measured lamination described by a measure $\mu\in{\mathbb R}^\Delta$.

\begin{figure}[h!]
\begin{center}
\epsffile{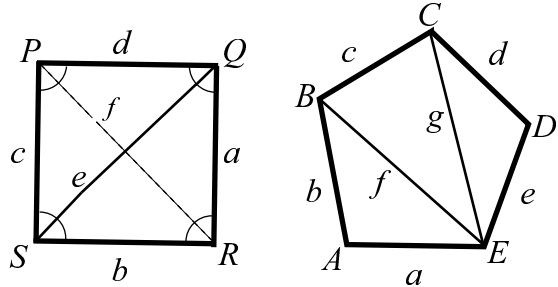}
\caption{Standard notation for squares and pentagons}
\label{fig:notation}
\end{center}
\end{figure}

One might hope to ``simplify'' the situation by choosing another ideal triangulation with smaller total intersection number $\sum_{X\in\Delta} \mu (X)$.  The example with three bands of arcs of respective weights 7,8,2 parallel to the corners $B,C,D$ of the pentagon in the notation of Figure \ref{fig:notation} (i.e., the respective measure on $a,\ldots ,g$ is given by 0,7,15,10,2,7,8)
shows that there are not always flips to decrease the total intersection number.  

A more delicate objective function to minimize is as follows.
Given a triangle $t$ with edges $a,b,c$, define
the {\it weight} or {\it perimeter} of $\mu$ on $t$ to be $\mu(t)={{\mu(a)+\mu (b) +\mu (c)}}$
and let
$$
\mu (\Delta)=\mu ({ T})=~{\rm sort}_\downarrow\{ \mu(t):t\in {T}\}
$$
be the sorted list of perimeters of triangles in ${T}$
ordered lexicographically.

\begin{lemma} \label{L:corner}
Consider
a flip along an edge $e$ of $\Delta$ triangulating a quadrilateral $\mathcal Q$ with corners $P,Q,R,S$, where e is incident on corners $P$ and $R$, as illustrated on the left in Figure \ref{fig:notation} and fix a measured lamination $\mu$ of $\mathcal Q$ transverse to the boundary.  The following conditions are equivalent:
\vskip 2mm
\noindent { i}) the flip along $e$ is non-decreasing for $\mu({ T})$;
\vskip 2mm

\noindent { ii}) $Q$ or $S$ has 
the least arcs around the corners of $\mathcal Q$, i.e., we have
$\aligned
\hskip 10mm min\biggl \{min\{a+e-b,e+d-c\}, &~min\{b+e-a,c+e-d\}\biggr\} \\
&\leq ~~~~min\biggl \{ a+b-e,c+d-e\biggr \};\\
\endaligned$

\vskip 2mm

\noindent { iii}) we have $ 2e+ max \{ a+b,c+d\} \leq max \{ a+c,b+d\} + max \{ a+d,b+c\}$;

\vskip 2mm 

\noindent { iv}) we have $e+max \{ a+b,c+d\} \leq f+max \{a+d,b+c\}$, i.e., the tropical face condition {\rm TF} holds.

\end{lemma}

\begin{proof}
(${ i} \Rightarrow { ii}$)~It is a small exercise with measures to confirm 
that the width of the band of leaves parallel to $P$ (or $R$) is ${{c+d-e}\over 2}$ (or ${{a+b-e}\over 2}$) and 
furthermore
the first (or second) term in the minimum on the left-hand side in part ${ ii})$ is twice the the width of the band of leaves which are parallel to $Q$ (or $S$).  Thus, the minimization in part ${ ii})$ is indeed equivalent to the condition on corners, which is
obviously a consequence of condition ${ i})$.

\vskip 2mm

\noindent (${ ii}\Rightarrow { iii}$)~If $min\{ a+e-b,e+d-c\}\leq min\{ b+e-a,c+e-d\}$, for example, then
$min\{ a+e-b,d+e-c\} \leq min\{ a+b-e,c+d-e\}$ implies, in particular, that $e\leq b$ and $e\leq c$.
Thus,\\
$$\aligned
&max \{ a+c,b+d\} + max \{ a+d,b+c\} \geq a+b+2c\geq 2e+a+b,\\
&max \{ a+c,b+d\} + max \{ a+d,b+c\} \geq 2b+c+d\geq 2e+c+d,\\
\endaligned$$\\
so\\
$$\aligned
max \{ a+c,b+d\} + max \{ a+d,b+c\}&\geq max \{ 2e+a+b,2e+c+d\}\\
&=2e+max \{ a+b,c+d\}\\
\endaligned$$
as required.  The remaining case that the band at $S$ is thinner than the band at $Q$ is analogous.

\vskip 2mm

\noindent (${ iii} \Leftrightarrow { iv}\Leftrightarrow { i}$)  Conditions ${ iii})$ and ${ iv})$ are actually equivalent by the tropical Ptolemy equation, and ${ iv)}$ is evidently equivalent to ${ i})$ by definition.
\end{proof}

It is worth emphasizing that the lemma shows in particular that the weight function gives the correct tropicalization TF of the convexity condition CF.

\begin{theorem}\label{T:seq}  Fix a measured lamination $\mu$ on  surface $F$.
If $\mu({ \Delta})$ is not minimal among all ideal  triangulations $\Delta$ of $F$, then there is a (possibly empty) 
finite sequence of weight-neutral flips followed by a flip that reduces it.
\end{theorem}

\noindent The proof will show that a global minimum exists in the special case that all leaves of $\mu$ are closed but not in general. Furthermore, one can bound above the number of neutral flips required before a reducing one can be found.

\begin{proof}
To begin, we assume that $\mu$ is generic in the sense that there are no coincidences among the values taken on the edges of $T$, i.e., no integral equalities hold. In this case, the theorem follows from the following

\begin{lemma}[Diamond Lemma]
Suppose that $T_1$ is a triangulation and each of $T_2,T_3$ arise from $T_1$ by a single weight-reducing flip. Then there is another triangulation $T_4$ which arises  from each of $T_2,T_3$ by a single weight reducing flip, i.e., the top
\begin{tabular}{ccccccc}
&&$T_1$\\
&$\swarrow$&&$\searrow$\\
$T_2$&&&&$T_3$\\\end{tabular}
implies the bottom
\begin{tabular}{ccccccc}
$T_2$&&&&$T_3$\\
&$\searrow$&&$\swarrow$\\
&&$T_4$\\
\end{tabular}
of a diamond of weight reducing flips.

\end{lemma}

\noindent First of all, this lemma is enough to prove the theorem.  Indeed, suppose that $\mu(T)>\mu(T')$  The triangulations $T$ and $T'$ can be connected by a sequence or path of flips.  If the flip starting from $T$ reduces the weight, then we are done.  Otherwise, there is some local maximum $T_1$ along the path of flips.  The diamond lemma allows us to push the 
weight of the path down near $T_1$.  Taking $T_1$ to be the closest maximum to $T$ along the path, then pushing down the weight at $T_1$ is tantamount to moving this first maximum closer to $T$.  This process can be repeated until the path is weight-decreasing immediately already at $T$ as required for Theorem \ref{T:seq}

\vskip 2mm

\noindent{\it Proof of the Diamond Lemma.} The argument is by cases on how the two flips from $T_1$ interact, and we consider the two quadrilateral neighborhoods supporting these two flips, comprised of at most four triangles.

\vskip 2mm

\noindent {\it Case 1:}~ If the two quadrilaterals are disjoint, then the two flips commute, and $T_4$ arises from $T_1$ by their composition.

\vskip 2mm

\begin{figure}[!h]
\begin{center}
\epsffile{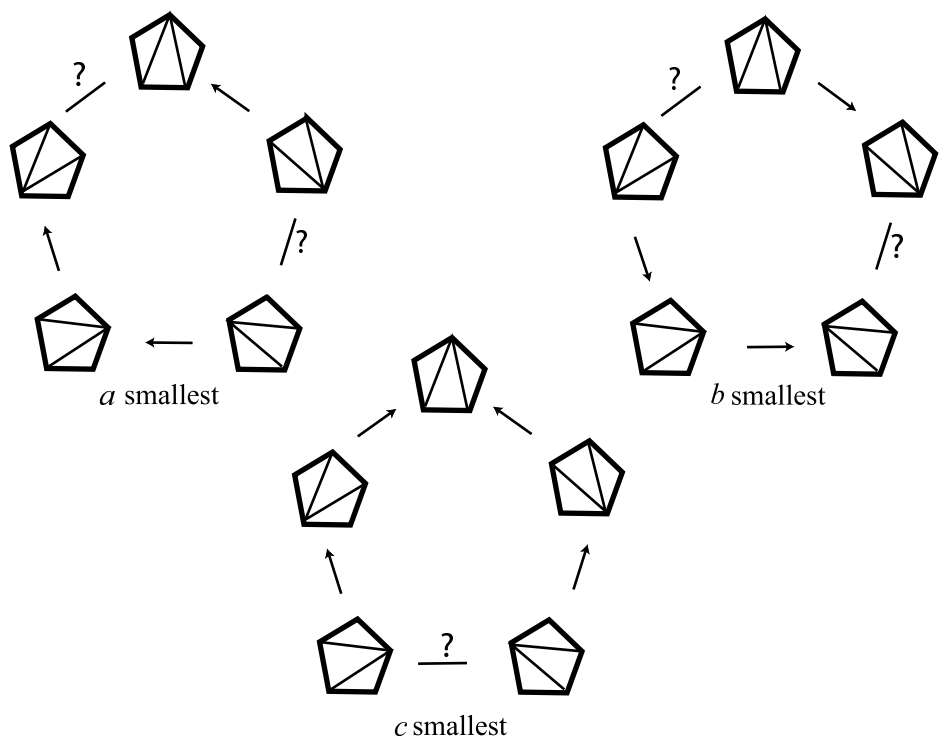}
\caption{Paths in a pentagon for Case 2}
\label{fig:case2}
\end{center}
\end{figure}

\noindent {\it Case 2:} If the two quadrilaterals overlap in one triangle, then their union is a pentagon with vertices labeled as illustrated in Figure \ref{fig:notation} but here letting $a,\cdots ,e$ denote the respective widths of bands of a measured foliation parallel to $A,\cdots ,E$.
We must check that regardless of the weights, there is a unique maximum and a unique minimum of weight among the five possible triangulations of the pentagon and all weight-nonincreasing flips, and we proceed by cases on which of $a,b,c,d,e$ is smallest.
 
\vskip 2mm

\noindent{\it a is smallest:}~  In this case, Lemma \ref{L:corner} forces the weight-decreasing flips illustrated in Figure \ref{fig:case2}a.
No matter how the remaining arrows are completed, there can be only one maximum or minimum.  A similar argument holds when $e$ is smallest.

\vskip 2mm

\noindent {\it b is smallest:}~ Now the forced weight-decreasing  flips are illustrated in Figure~\ref{fig:case2}b, and again, there are unique extrema.  A similar argument holds when $d$ is smallest.

\vskip 2mm

\noindent {\it c is smallest:}~ In this case, almost everything is fixed as illustrated in Figure~\ref{fig:case2}c.

\vskip 2mm

\noindent {\it Case 3:}~ Suppose the two quadrilaterals overlap in two triangles configured as an annulus
as illustrated in Figure \ref{fig:case3}. 
We may assume without loss of generality that the measured lamination is in its standard position as illustrated
on the top-right in Figure~\ref{fig:case3}.  The triangulations are then necessarily ordered in the natural way
as illustrated on the bottom in the figure.  Thus, there can be no local maximum, and this contradicts the initial assumptions.

\vskip 2mm

\begin{figure}[h!]
\begin{center}
\epsffile{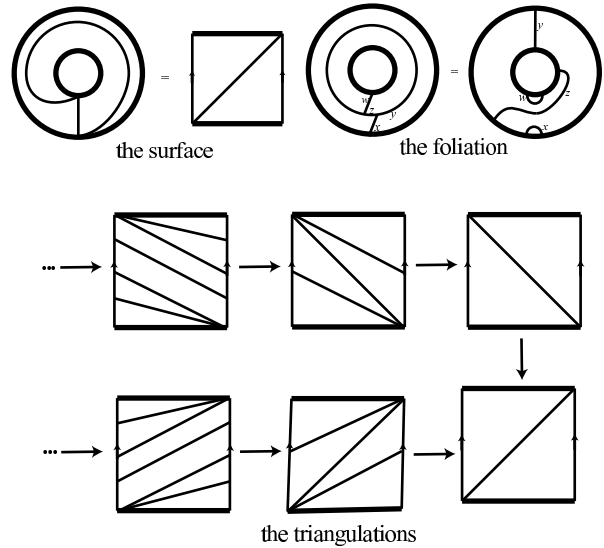}
\caption{No local maximum for an annulus in Case 3}
\label{fig:case3}
\end{center}
\end{figure}

\vskip 2mm

\noindent {\it Case 4:}~ Suppose the two quadrilaterals overlap in two triangles, which is configured as a bigon
as illustrated in Figure \ref{fig:case4}.
  Again, we may assume without loss of generality that the measured lamination is standard as illustrated on the top-right in Figure \ref{fig:case4}.  On the bottom of this figure are shown the two possible transitions from the unique possible initial triangulation each of whose edges has distinct endpoints, and arithmetic shows that the two transitions are weight-decerasing under the conditions depicted.  However, it cannot be that both $x-y$ and $y-x$ are negative, so the initial triangulation
is not a local maximum again contradicting the assumptions.

\vskip 2mm

\begin{figure}[h!]
\begin{center}
\epsffile{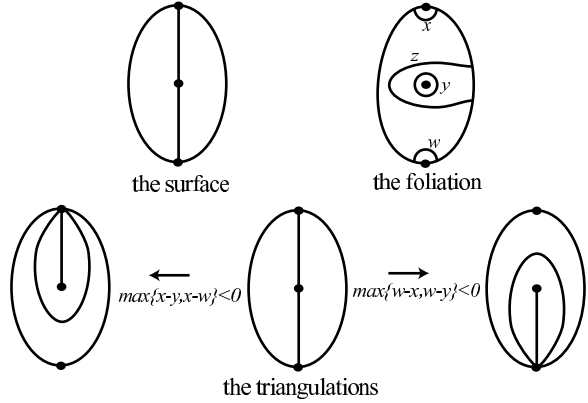}
\caption{No local maximum for a bigon in Case 4}
\label{fig:case4}
\end{center}
\end{figure}

This completes the proof provided the lamination is generic.
In the general case, say $\mu$ is an integral measured lamination and let $\mu'=\mu+\varepsilon$ be a generic perturbation of $\mu$.  We may perform flips trying to improve with respect to $\mu'$ as much as possible.  These flips will either be weight-decreasing or -neutral.  However, there is a bound in terms of the genus on how many neutral flips can be performed in sequence before returning to the same combinatorial type of triangulation, and hence there is a mapping class that leaves $\mu$ invariant.  Such a sequence can therefore be removed from an overall weight-decreasing sequence of flips, so the general case follows from the generic one.
\end{proof}

Together with part c of Lemma \ref{implications}, this theorem has the following immediate consequence:

\begin{corollary}\label{c:exist}
For any decorated measured lamination $\mu$ with all leaves closed in a punctured surface, there is a (not necessarily unique) ideal triangulation $\Delta$ with complementary triangles $T$ so that $\mu$ satisfies ${\rm TT}$ on $\Delta$.  Moreover, there is  an algorithm for finding $\Delta$
based on minimizing first the sorted perimeters 
${\rm sort}_\downarrow\{ \mu(a)+\mu(b)+\mu (c):t\in T\}$
of triangles and then the total weight $\sum_{e\in\Delta}\mu (e)$ of edges.
\end{corollary}

\bibliographystyle{amsplain}

 \end{document}